\documentclass[colorinlistoftodos,a4paper,a4wide,11pt]{article}
\def\one{1}

\def\two{2}

\def\hidereviews{2}

\usepackage{xcolor}
\usepackage{listings}
\definecolor{eclipseStrings}{RGB}{42,0.0,255}
\definecolor{eclipseKeywords}{RGB}{127,0,85}
\colorlet{numb}{magenta!60!black}

\lstdefinelanguage{json}{
  basicstyle=\normalfont\ttfamily,
  commentstyle=\color{eclipseStrings}, %
  stringstyle=\color{eclipseKeywords}, %
  numbers=left,
  numberstyle=\scriptsize,
  stepnumber=1,
  numbersep=8pt,
  showstringspaces=false,
  breaklines=true,
  frame=lines,
  backgroundcolor=\color{white}, %
  string=[s]{"}{"},
  comment=[l]{:\ "},
  morecomment=[l]{:"},
  literate=
  *{0}{{{\color{numb}0}}}{1}
  {1}{{{\color{numb}1}}}{1}
  {2}{{{\color{numb}2}}}{1}
  {3}{{{\color{numb}3}}}{1}
  {4}{{{\color{numb}4}}}{1}
  {5}{{{\color{numb}5}}}{1}
  {6}{{{\color{numb}6}}}{1}
  {7}{{{\color{numb}7}}}{1}
  {8}{{{\color{numb}8}}}{1}
  {9}{{{\color{numb}9}}}{1}
}

\usepackage{tikz}
\usepackage{pgfplots}
\pgfplotsset{compat=1.14}

\usepackage{geometry}
\usepackage[utf8]{inputenc}
\usepackage[T1]{fontenc}
\usepackage{amsmath,amssymb,commath,mathtools,amsfonts}
\usepackage{mathrsfs}
\usepackage{enumitem}
\usepackage{booktabs}
\usepackage[margin=0.1\textwidth]{caption}
\usepackage{siunitx}
\usepackage[normalem]{ulem}
\usepackage{todonotes}

\newcommand{\review}[3][{noinline}]{\if\hidereviews\one#2\else%
  \if\hidereviews\two{\color{red}#2}%
  \else{\color{red}#2}\todo[#1]{REVIEW}\sout{#3}\fi%
  \fi}

\usepackage[pdfusetitle]{hyperref}

 \hypersetup{
   colorlinks=true,
  citecolor=blue,
  linktocpage=true,
  pdfstartview=FitV,
  breaklinks=true,
  pdfpagemode=UseNone,
  pageanchor=true,
  pdfpagemode=UseOutlines,
  plainpages=false,
  bookmarksnumbered,
  bookmarksopen=true,
  bookmarksopenlevel=1,
  hypertexnames=false,
  pdfhighlight=/O,
  nesting=true,
  pdfsubject={Mathematics, Numerical analysis},
  pdfkeywords={power networks, gas networks, numerical simulation},
  pdfcreator={pdfLaTeX},
  pdfproducer={},
}

\usepackage{amsthm,thmtools}

\usepackage[backend=bibtex,style=alphabetic,sortcites=true,giveninits=true]{biblatex}  %
\DeclareMathOperator{\Rey}{Re}
\DeclareMathOperator{\ld}{ld}

\let\abs\relax
\DeclarePairedDelimiter\abs{\lvert}{\rvert}

\newcommand{\im}{\text{i}}

\title{Efficient simulation of coupled gas and power networks under uncertain demands}
\author{Eike Fokken\footnotemark[1] \and Simone G\"ottlich\footnotemark[1] \and Michael Herty\footnotemark[2]}

\begin{document}

\maketitle
\renewcommand*{\thefootnote}{\fnsymbol{footnote}}
\footnotetext[1]{Department of Mathematics, University of Mannheim, 68131 Mannheim, GERMANY, \{\href{mailto:fokken@uni-mannheim.de}{fokken},\href{mailto:goettlich@uni-mannheim.de}{goettlich}\}@uni-mannheim.de}
\footnotetext[2]{IGPM, RWTH Aachen University, Templergraben 55. 52056 Aachen, GERMANY, \href{mailto:herty@igpm.rwth-aachen.de}{herty@igpm.rwth-aachen.de}}

\renewcommand*{\thefootnote}{\arabic{footnote}}%

\renewcommand*{\sectionautorefname}{Section}

\let\subsectionautorefname\sectionautorefname
\let\subsubsectionautorefname\sectionautorefname

\abstract{We introduce an approach and a software tool for solving coupled energy networks composed of gas and electric power networks. Those networks are coupled to stochastic fluctuations to address possibly fluctuating demand due to fluctuating demands and supplies. Through computational results the presented approach is tested on networks of realistic size.

}

\section{Introduction}
\label{sec:introduction}

In view of the changing energy demands and supplies the combined and intertwined energy networks are expected to play a prominent role in the future. Here, the additional but unpredictable and volatile  energy sources need to be complemented with traditional means of production as well as possibly additional large--scale storage~\cite{IEA2019,EESI2019}. Following previous approaches~\cite{CHERTKOV2015541,ZengFangLiChen2016,Zlotnik16a} we consider here energy storage through coupling power networks to gas networks. The latter are able to generate sufficient power at times when renewable energy might not be available or vice versa convert energy to ramp up pressure in the gas networks, see e.g. \cite{Schiebahn15,Brown18a,Heinen16}. A major concern when coupling gas and power networks is guaranteeing a stable operation even at times of stress due to (uncertain) heavy loads. The propagation of 
possible uncertain loads on the power network and its effect on the gas network has been subject to recent 
investigation and we refer to \cite{CHERTKOV2015541} and references therein. Contrary to the cited reference \cite{CHERTKOV2015541} we are interested here in a full simulation of both the gas and the power network as well as the a simulation of the stochastic demand, respectively supply. This will allow for a prediction at all nodes as well as study dynamic  effects changing supplies and demands in the network.  
\par 
Regarding the underlying models for simulation we rely on established power flow, gas flow and stochastic demand models that will be briefly reviewed in the following. Power flow (PF) is typically modeled through prescribing real and reactive power at nodes of the electric grid. Their values are obtained through a nonlinear system of algebraic equations. Supply and demand can be time--dependent requiring to 
frequently resolve the nonlinear system. For more details on the model we refer to \cite{Frank16,Capitanescu16a,kit:faulwasser18d,Low14a,Low14b,kit:muehlpfordt18d,Bienstock14a} as well as
to the forthcoming section where the equations are reviewed. While propagation of electricity is typically assumed to be instantaneous as in the PF equations, the propagation of gas in networks has an intrinsic spatial and temporal scale. A variety of models exist nowadays and we follow here an approach based on hyperbolic balance laws as proposed e.g. in \cite{coupling_conditions_isothermal,banda_herty_klar,bressan_canic_garavello_herty_piccoli,brouwer_gasser_herty,doi:10.1137/060665841}. This description allows the prediction of gas pressure and gas flux at each point in the pipe as well as nodes of the network. Both quantities are relevant to assess possible stability issues as well as allow for coupling towards the electricity network. The numerical solution of the governing gas equations as well as their coupling towards electricity and towards other gas pipelines is also detailed in the forthcoming section. Finally, we recall recent results on modeling of the prediction of the electricity demand that will be used to simulate the uncertain power fluctuations. Here, we follow models introduced in 
\cite{Aid.2009, Kiesel.2009, LuciaSchwartz.2002, SchwartzSmith.2000, Wagner.2014,Barlow.2002} and  the monograph \cite{Benth.2008} that prescribe the electricity demand as Ornstein--Uhlenbeck processes. Let us emphasize that our approach is not limited to the particular application of gas transportation but could eventually be applied to problems of traffic flow and supply chain dynamics on networks.  
\par 
Finally, we point to other existing numerical simulations with possibly similar objective. For example, in \cite{DissKolb} a solver for a gas network based on (general) hyperbolic balance laws has been introduced. This implementation serves as foundation for the method introduced in the forthcoming section. Our tool has the advantage of easier extensibility, an open source license and a modern design approach featuring for example extensive software testing.  In addition we also include stochastic power demands in the powerflow network setting. In \cite{AndersonTurnerKoch2020} uncertainty in PF is computed relying on approaches based on neural networks. A difference to the presented approach is the restriction to linear powerflow problems and the absence of coupling to gas networks.
A further concept plan4res, see \cite{Beulertz:774907} has been presented to also address general renewable energy sources as well as energy distribution based on discrete optimization approaches, which focuses on energy system modeling in more generality.
Furthermore, there exists a software suite by Fraunhofer SCAI called MYNTS, see \cite{simultech16}, that also includes simulation and optimization of inter--connected grid operations with a focus on the design of suitable networks.  %

The paper is organized as follows: In \autoref{sec:model} we introduce the model equations for the coupled network setting and the uncertain demand. The numerical discretization is then given in \autoref{sec:discretization}. \autoref{sec:scenario} is concerned with the presentation of our software suite and the numerical investigation of relevant scenarios.

\section{Mathematical modeling}%
\begin{figure}[ht]
  \centering
  \includegraphics[width=0.9\textwidth]{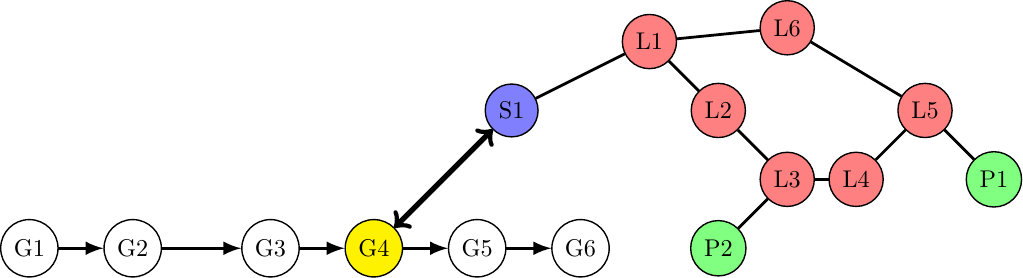}
  \caption{A schematic example of the kind of network under consideration.
    The upper right part is a power network with blue slack nodes, green powerplants and red load nodes.
    In the lower left there is a gas network with pipelines between junctions.
  The doubly-pointed arc is a gas-power connection.}
  \label{fig:networkschema}
\end{figure}

\label{sec:model} 
In this section we introduce the mathematical models to be discretized in the forthcoming section. 
We denote by $\mathcal{G}=(\mathcal{N},\mathcal{A})$ a directed graph with a set of nodes $\mathcal{N}$ and a set of arcs $\mathcal{A}$.
All dynamics are either on nodes and/or on arcs of the graph.
The whole graph is sub-divided into a power network, a gas network and a set of arcs connecting these two:
\begin{equation*}
  \begin{aligned}
    \mathcal{G} &= \mathcal{G}_\text{P} \cup\mathcal{G}_\text{G} \cup\mathcal{G}_\text{GP}\\
    \mathcal{N} &= \mathcal{N}_\text{P} \cup\mathcal{N}_\text{G} \cup\mathcal{N}_\text{GP}\\
    \mathcal{A} &= \mathcal{A}_\text{P} \cup\mathcal{A}_\text{G} \cup\mathcal{A}_\text{GP}.\\
  \end{aligned}
\end{equation*}
The different parts behave quite differently.  In the power network the arcs just carry two parameters and their topological information, i.e.~their starting and ending node and the nodes carry most of the physical information, namely active and reactive power, while the situation in the gas network is reversed.
Here the arcs carry a balance law describing gas dynamics while the nodes only carry coupling information.
Yet in all parts of the network only the nodes have (possibly stochastic) boundary conditions\footnote{In the power network these are rightfully called \emph{node specifications} and we only call them boundary conditions to unify the wording between gas and power networks.}.
The gas-power connection part of the network consists only of arcs which relate power demand and gas consumption or gas generation and power surplus.
An illustration of such a network can be found in \autoref{fig:networkschema}.

Next, we define model equations for each node and arc (where applicable) based on physical models. 

\subsection{Power flow equations modeling electric power flow on $\mathcal{G}_P$}%
\label{sec:power-model}
The evolution of reactive and active power at nodes is modeled by the power flow equations \cite{grainger2016power}. This model can be used to describe the behavior of power networks operating
at sinusoidal alternating current (AC)~\cite{fokken2021relation}. The quantities modeled are the active or real power $P_k=P_k(t)$ and the reactive power $Q_k=Q_k(t)$ present at each node $k\in\mathcal{N}_\text{P}$ at time $t$.
Those  are  functions of the voltage magnitude $V_k$ and angle $\phi_k$.
Further,  we model the  admittance of each component, denoted
by $Y$, which is written  into real and imaginary part \(Y = G+\im B\).
The admittance is the inverse of the impedance which in turn is a complex extension of Ohmic resistance in the power network.

The admittance of a transmission line, that is, an arc \(a \in \mathcal{A}_\text{P}\) connecting nodes $i,k \in \mathcal{N}_P$
is denoted by \(Y_{ik} = G_{ik}+\im B_{ik}\), which we set to zero, if no arc connects \(i\) and \(k\).
The admittance of a node \(k \in \mathcal{N}_P\) is denoted by \(Y_{kk} = G_{kk}+\im B_{kk}\).

With this, we can write down the power flow equations, a set of \(2 | \mathcal{N}_P | \)
equations at any point $t$ in time of the type
\begin{equation}
  \label{eq:1}
  \begin{aligned} P_{k}(t) & = \sum_{i \in \mathcal{N}_\text{P}}
    V_k(t)V_i(t)(G_{ki}\cos(\phi_k(t)-\phi_i(t))+B_{ki}\sin(\phi_k(t)-\phi_i(t))),\\
    Q_{k}(t) & = \sum_{i \in \mathcal{N}_\text{P}}
    V_k(t)V_i(t)(G_{ki}\sin(\phi_k(t)-\phi_i(t))-B_{ki}\cos(\phi_k(t)-\phi_i(t))),
  \end{aligned}
\end{equation}
for the unknowns $(P_k(t),Q_k(t),V_k(t),\phi_k(t)) k \in \mathcal{N}_P$.  In order to obtain a unique solution, additional $2| \mathcal{N}_P |$ equality constraints have to be specified.  We distinguish three different equality constraints according to the type of node $k\in \mathcal{N}_P$:
\begin{itemize}
\item Slack nodes $k$ specify values for the voltage magnitudes and angles \(V_k,\phi_k\).
\item Load nodes  $k$ specify values for the active and reactive power \(P_k,Q_k\).
\item Generators  $k$ specify values for the active power and voltage magnitude \(P_k,V_k\).
\end{itemize}
A necessary condition for uniqueness of the power flow equations is to have at least one slack node,
because otherwise for any solution \((P_k(t),Q_k(t),V_k(t),\phi_k(t)) k \in \mathcal{N}_P\) of the power flow equations, a second one is obtained by shifting all phase angles \(\phi_k\) by the same amount \(r\in\mathbb{R}\).  This is possible, as without a slack node, the equations only depend on differences in the phase angles.
Often only a single slack node is used, although also multiple slack nodes can be used \cite{Chayakulkheeree2007}.

Instead of the described AC powerflow equations, it is possible to use so-called direct current (DC) powerflow equations, which are a linear approximation, see \cite[6.10]{10.5555/559961} for an overview.
This approach simplifies the numerical treatment greatly at the cost of some accuracy.
Other linearizations are subject of active research, see e.g.~\cite{LiPanLiu2019}.

\subsection{Mathematical modeling of gas flow on $\mathcal{G}_\text{G}$}%
\label{sec:gas-model}
We model the following quantities of the gas flow, namely, the pressure $p=p(t,x)$ as well as the flux $q=q(t,x)$. The units of those quantities are (\(\si{\bar}\)) and  \(\si{\meter \cubed \per \second}\).  Note that we use the volumetric flow as opposed to mass flow, as is customary in real-world gas networks.
The pressure is given by a function of the gas density $\rho=\rho(t,x)$. An overview as well as recent results on gas flow can be found e.g. in \cite{coupling_conditions_isothermal,banda_herty_klar,bressan_canic_garavello_herty_piccoli,brouwer_gasser_herty,doi:10.1137/060665841,reigstad,reigstad2}

\subsubsection{Transport of gas along pipelines}%
\label{sec:edges}
The direction of the arcs determines the positive direction of the flow.
We distinguish two types of arcs, pipelines and controlled arcs, which are described in the next \autoref{sec:controlled-gas-arcs}.
In contrast to the power arcs, pipelines \(a \in \mathcal{A}_\text{G}\) in the gas
network are modeled as an interval \([0,L_a]\) with further structure.
The gas flow in pipelines is modeled 
with the {isentropic Euler equations} as e.g. proposed in
\cite{banda_herty_klar}.

The conservative variables are the gas density \(\rho\) and the flux $q$.
Consider an arc $a \in \mathcal{A}_\text{G}$ parametrized by $x \in (0,L_a)$.
Then, the density $\rho=\rho_a(t,x)$ and $q=q_a(t,x)$ fulfill in the weak sense the
following system of hyperbolic balance laws for each $t\geq 0$.
\begin{equation}\label{eq:gas}
  \begin{pmatrix}\rho\\  q\end{pmatrix}_t + \begin{pmatrix} \frac{\rho_0}{A}q\\ \frac{A}{\rho_0}p(\rho) + \frac{\rho_0}{A}\frac{ q^2}{\rho}\end{pmatrix}_x = \begin{pmatrix}0\\S(\rho, q) \end{pmatrix}.
\end{equation}
Here \(S\) is the source term modeling wall friction in the pipes detailed
below, \(A\) is the cross-section of the pipe, \(\rho_0\) is the (constant)
density of the gas under standard conditions and \(p\) is the pressure function, which we describe in detail below.
First we note that the system \eqref{eq:gas} is accompanied by initial conditions 
\begin{align}
	\rho(0,x)=\rho_{0,a}(x), \; q(0,x)=q_{0,a}(x)
\end{align}
that may also dependent on the selected arc $a$ and describe the initial state of density and flux.
Suitable boundary conditions at $x=0$ and $x=L_a$ will be discussed below when the coupling at nodes of the network is introduced.

The pressure  \(p\) is a function of the density given by
\begin{equation}
\label{eq:pressure}
  p(\rho) = \frac{c_{\text{vac}}^2\rho}{1-\alpha c_{\text{vac}}^2 \rho},
\end{equation}
where \(c_{\text{vac}}\) is the vacuum limit (\(\rho \to 0\)) of the speed of sound and
\(\alpha\) is a measure of compressibility of the gas. The relevant constants of
the gas network are gathered in \autoref{tab:standardconditions}.
It is possible to express the density in terms of the pressure as
\begin{equation}
  \label{eq:density}
  \rho = \frac{p}{c_{\text{vac}}^2 z(p)},
\end{equation}
where \(z(p) = 1+\alpha p\) is the so-called compressibility factor.
\begin{table}[ht]
  \centering
  \caption{Gas net constants.}
  \begin{tabular}{llllll}
  	\toprule
    \(\rho_0 \, [\si{\kilogram\per\metre\cubed}]\) &  \(c_\text{vac}\, [\si{\meter \per \second}]\) &\(\alpha \, [\si{\per \bar}]\)\\
    \midrule
    0.785 &   364.87  &     -0.00224 \\
    \bottomrule
  \end{tabular}%
  \label{tab:standardconditions}
\end{table}

The source term \(S\) is given as
\begin{equation}
  \label{eq:13}
  S(\rho,q) = \frac{\lambda(q)}{2d} \frac{\abs{q}}{\rho} \left( -q \right),
\end{equation}
where now \(d\) is the pipe diameter and \(\lambda(q)\) is the flux-dependent Darcy friction factor, see \cite{Brown_history_darcy_weisbach}.
The friction is governed by the so-called Reynolds number,
\begin{equation}
  \label{eq:15}
  \Rey(q) = \frac{d}{A\eta}\rho_0\abs{q}
\end{equation}
(here \(\eta = \SI{e-5}{\kilogram\per \meter \per \second}\) is the dynamic viscosity of the gas).
For \(\Rey <2000 \) the friction is dominated by laminar flow and according to \cite{Menon2015} we may assume
\begin{equation}
  \label{eq:16}
  \lambda(q) = \frac{64}{\Rey(q)}.
\end{equation}
For \(\Rey > 4000 \) the friction is dominated by turbulent flow, see again \cite{Menon2015}, and the Swamee-Jain
approximation~\cite{swamee-jain} is used, i.e., 
\begin{equation}
  \label{eq:17}
  \lambda(q) = \frac{1}{4} \frac{1}{\ld(\frac{k}{3.7 d} + \frac{5.74}{\Rey^{0.9}})^2}.
\end{equation}
For $\Rey$ in the intermediate regime the numbers are interpolated using a cubic polynomial differentiable at $\Rey \in \{ 2000, 4000 \}$.

\subsubsection{Controlled gas arcs}
\label{sec:controlled-gas-arcs}
In addition to pipes, the considered gas network also contains a compressor and a control valve.
Both are modeled similarly and come equipped with a control function \(u\) that influences the pressure.

Both valves and compressors do not influence the flow rate of the gas inside them, so at every timepoint \(t\) there must hold
\begin{equation}
  \label{eq:18}
  q_\text{out}(t) = q_\text{in}(t).
\end{equation}
Yet for the pressure we set in compressors
\begin{equation}
  \label{eq:19}
  p_\text{out}(t) = p_\text{in}(t)+u(t),
\end{equation}
and in valves
\begin{equation}
  \label{eq:20}
  p_\text{out}(t) = p_\text{in}(t)-u(t),
\end{equation}
such that the control can change the pressure. For other possible compressor models see for example \cite{DissKolb}.  In addition we demand \(u(t)\geq 0\) for both component types. Also for the purpose of optimization, using the compressor comes with a cost, while using the valve is free.
Note that for the main part of this work, both valves and compressors are in-active meaning a control function of \(u(t)=0\).
Only in the optimization example a non-zero control is allowed.

\subsubsection{Nodes of the gas network $\mathcal{N}_G$}%
\label{sec:nodes}
The previous set of differential equations has to be accompanied by boundary conditions if $L_a<\infty.$ 
For nodes $n \in \mathcal{N}_G$ coupling only gas pipelines those are typically described in terms of coupling conditions \cite{Herty2007}. Those yield an implicit description of the boundary values in terms of physical relations. Several different conditions exist, see \cite{reigstad,reigstad2,holle2020new}.
Yet, for the purpose of real-world gas pipelines,
\cite{Fokken2021} seems to indicate that coupling via equality of pressure is sufficient
for the expectable accuracy of the whole modeling approach.
Therefore consider a node \(n \in \mathcal{N}_G \) with \(K\) adjacent arcs.
Let $E  \subset \mathcal{A}$ denote the set of adjacent arcs and let
\begin{equation}
  \label{eq:3}
  \begin{aligned}
    &s: E \to \set{\pm 1},\\
    &s(e) =
    \begin{cases}
      \phantom{-{}}1\ &e \text{ starts in } n\\
     -1\ &e \text{ ends in } n\\
    \end{cases}
  \end{aligned}
\end{equation}
distinguish arcs starting and ending in \(n\).
Also let \(p_e(t),q_e(t)\ e \in E\) denote the boundary values at time $t$ of arc \(e\) in node \(n\), i.e.,
$p_e = p_e(t,x=0)$ if $s(e)=1$ and $p_e(t)=p_e(L_a,t)$ if $s(e)=-1$ and similarly for $q_e.$ Then, the coupling and boundary conditions at node \(n\) read
\begin{equation}
  \label{eq:4}
  \begin{aligned}
    p_e(t) &= p_f(t) \text{ for all } e,f \in E \\
    q_n(t) &= \sum_{e \in E}s(e)q_e(t),
  \end{aligned}
\end{equation}
where \(q_n:\mathbb{R}^+_0 \to \mathbb{R}\) a possibly time dependent external and given demand or supply function.  These are in total \(\abs{E}\) equations at node \(n\) at each point in time.

For a single node with $K$ adjacent arcs extending to infinity and under
subsonic condition for the initial data existence of weak entropic solutions in
$BV$ has been shown e.g. in \cite{Colombo2008}. In \cite{holle2020new} existence
of weak integrable solutions on a graph has been established. Similar results
are also available for other choices of coupling conditions and we refer to e.g.
\cite{BressanCanicGaravelloHertyPiccoli2014}.

\subsection{Modeling of Gas-to-Power or Power-to-Gas nodes $\mathcal{N}_\text{GP}$}%
\label{sec:gas-power-conversion}
Gas power plants are also modeled as arcs in the graph connecting a power node
and a gas node. They transform gas into power at a linear rate and also power
into gas at a (different) linear rate as was done in~\cite{Fokken2021}. At the switching point we smooth the
resulting kink with a polynomial. This has no physical counterpart and is done purely for numerical reasons.
Note that due to technical reasons our
polynomial maps gas flow to power instead of the other way around as was chosen
in~\cite{Fokken2021}.
For the power output of the gas power plant there holds
\begin{equation}
  \label{eq:5}
  P =
  \begin{cases}
    E_{\text{PtG}}\cdot q &\text{ for } \phantom{ {}- \kappa < {}} q<-\kappa\\
    \text{pol}(q) &\text{ for } -\kappa <q<\kappa \\
    E_{\text{GtP}}\cdot q &\text{ for } \phantom{{}-{}}\kappa<q 
  \end{cases}\ ,
\end{equation}
where \(E_{\text{PtG}}\) and \(E_{\text{GtP}}\) are the efficiencies of the
power-to-gas process and the gas burning respectively and \(\text{pol}\) is the
interpolating polynomial. For \(\kappa\) we choose \(\kappa = \SI{60}{\meter
  \cubed \per \second}\). 
The power is taken as the real power of the attached power node, which is a slack node and therefore
provides via the solution of the powerflow problem a real power demand.

\subsection{Stochastic power nodes $\mathcal{N}_\text{P}$}%
\label{sec:stoch-power-demand}
In order to incorporate uncertain power demands into our model we add a new kind
of load node, the \emph{stochastic PQ-node}. The type of uncertainty employed,
the Ornstein-Uhlenbeck process, has a long history in modeling uncertain demands
of various types and has also been used for electricity demand, see
\cite{Barlow.2002,GoettlichKoblLux2020}.
In \cite{GoettlichKoblLux2020} a setting similar to ours was examined
but applied to the Telegrapher's equations instead of the powerflow equations.

The stochastic PQ-node, just like its deterministic cousin prescribes a real and
reactive power demand as boundary conditions but now these demands are
stochastic time-dependent quantities modeling the uncertainty of demand at this
node. Of course this uncertainty is not total, as one may expect the demand to
follow historic timelines of demand or some other estimate derived from
knowledge about the season, weather or even current events like a sports
tournament. This structure of uncertain fluctuation about a deterministic
estimate suggests using a mean reverting stochastic process for the power demand,
\(\left(P_t\right)_{t \in [0,T]}\), that is, a process that is drawn back to
some deterministic function \(\mu(t)\) over time. If we further assume that
fluctuations around \(\mu\) are independent of the current time and also of the
current value of \(P\), a natural choice for the process is the
Ornstein-Uhlenbeck process (OUP). It is characterized by the following
stochastic differential equation,
\begin{align}
	dP_t=\theta\left(\mu(t)-P_t\right)dt + \sigma dW_t,\ \quad P_{t_0}=p_{0}, \label{eq:oup}
\end{align}
where \(W_t\) is a one-dimensional Brownian motion, \(\theta,\sigma >0\) are the so-called \emph{drift} and \emph{diffusion} coefficients, and \(p_{0}\) is the demand at the starting time \(t_0\).

Whenever the current demand \(P_t\) differs from \(\mu(t)\), the drift term enacts a force towards the deterministic demand estimate \(\mu(t)\).
This behavior is called \emph{mean reversion}.
The size of this force is characterized by the drift coefficient \(\theta\).
In absence of diffusion (for \(\sigma = 0\)), the OUP degenerates to a deterministic ordinary differential equation, that is drawn to the mean exponentially.
For \(\sigma >0\) on the other hand, this mean reversion is counteracted by fluctuations, whose size is determined by \(\sigma\).
For images of OUP realizations see \autoref{fig:cutoff-oup}.

Both the real power demand \(P\) as well as the reactive power demand \(Q\) are realized as an OUP in our setting.

Note that it is even possible to solve the stochastic differential equation \eqref{eq:oup} explicitly via 
\begin{equation*}
	P_t = p_{0} e^{-\theta (t-t_0)} + \theta \int_{t_0}^te^{-\theta(t-s)} \mu(s)ds + \sigma \int_{t_0}^te^{-\theta(t-s)}dW_s.
\end{equation*}
From this explicit expression one can see that \(P_t\) is normally distributed with mean
\begin{equation*}
  \mu_t = p_{0}e^{ - \theta (t-t_0)}  + \theta \int_{t_0}^t e^{ - \theta \left( {t - s} \right)} \mu \left( s \right)ds
\end{equation*}
and variance
\begin{equation*}
  V = \sigma ^2 \int\limits_{t_0}^t {e^{ - 2\theta \left( {t - s} \right)}ds}\ .
\end{equation*}
The mathematical properties as well as the possibility to account for forecasts make the OUP a prime candidate for modeling uncertainty in power demand, see also~\cite{Benth.2008}.

\section{Discretization}%
\label{sec:discretization}
Having defined our model we now need to discretize it in order to search
solutions numerically.
This search will be carried out by Newton's method at each time step.
The discretization is different for pipelines (with their
balance law) and stochastic PQnodes on the one hand and all other components on
the other hand. This is because only pipelines and stochastic PQnodes couple the
state of the model at different times, because only they contain time derivatives.

Therefore we choose a time discretization with uniform stepsize \(\Delta t\) for
a time horizon \([t_\text{start},t_\text{end}]\) such that
\(J = \frac{t_\text{end}-t_\text{start}}{\Delta t}\) is an integer and henceforth consider only the discretized time points
\(j \Delta t\), where \(0\leq j \leq J\).
For all equations except the isentropic Euler equations and the
Ornstein-Uhlenbeck process, this means we simply evaluate them at the time
steps.

\subsection{Power discretization}%
\label{sec:power-discretization}
In the power network this means we evaluate equations~\eqref{eq:1} at each time step \(j\Delta t\).
Note that therefore we must also evaluate the boundary conditions at each time step for each node.

\subsection{Gas pipeline discretization}%
\label{sec:gas-discretization}
For the isentropic Euler equations we need a suitable numerical scheme.
For a pipeline of length \({L_a}\)
we introduce a space discretization with stepsize \(\Delta x_a\), such that
\(K \coloneqq \frac{{L_a}}{\Delta x_a}\) is an integer. We replace the
continuous values of pressure (or density, see equations~\eqref{eq:pressure},
\eqref{eq:density}) and flow with values at each \(x = x_k \coloneqq k\Delta
x_l\), \(0\leq k\leq K\). The isentropic Euler equations themselves are
discretized with an implicit Box scheme due to Bales
et.~al.~\cite{KolbLangBales2010}. For a general hyperbolic balance law
\begin{equation}
  \label{eq:2}
  u_t +{f(u)}_{x} = g(u)
\end{equation}
with space discretization \(x_k,\Delta x_a\) as above we have for the time
step between \(t\) and \(t^{*} = t+\Delta t\)
\begin{equation}
  \label{eq:8}
  \frac{u_{k}^{*}+u_{k-1}^{*}}{2} =  \frac{u_{k}+u_{k-1}}{2}
  -\frac{\Delta t}{\Delta x_a}\left(f(u_{k}^{*})-f(u_{k-1}^{*}) \right)
  + \Delta t\left( g(u_{k}^{*})+g(u_{k-1}^{*}) \right),
\end{equation}
where \(u_{k} = u(x_k,t)\) and \(u_{k}^{*} = u(x_k,t^{*})\). In our case
\(u_{k}\) has two components, density and flux and hence we get \(2K\)
equations on a pipeline for \(2K+2\) variables. Therefore for each pipeline we
need an additional \(2\) equations for the possibility of a unique solution,
namely one boundary condition for \(u_0\) and one for \(u_K\).

Note that no diagonalization is needed before a time step and equation \eqref{eq:8} can be used directly.
But an inverse CFL condition
\begin{equation*}
  \Delta t > \frac{\Delta x}{2\Lambda},
\end{equation*}
where \(\Lambda = \min\{\abs{\lambda(u)}\ |\lambda \text{ is an eigenvalue of
}f^\prime(u)\}\), must be fulfilled, which also shows that the scheme breaks down
for transonic flow, where an eigenvalue approaches \(0\). We refer to \cite[Prop 4.2,
following remark]{DissKolb} for a proof in the scalar case and \cite[section 4.1]{Fokken2019b}
for a numerical study of systems of conservation laws.

The inverse CFL condition is well-suited for the task at hand, as large time
steps are desirable for numerical feasibility when simulating over many hours.

These are of course supplied by the nodes, which yield a single equation for
each arc connected to them\footnote{Note that this situation breaks down
  whenever the flow in pipes is supersonic, see
  e.g.~\cite{MartinHertyMueller2017}.} 
from its starting node and one from its ending node.

We also remark that discretizing the controlled gas arc equations \eqref{eq:18}, \eqref{eq:19} and \eqref{eq:20} is straightforward.

\subsection{Node discretization}
\label{sec:node-discretization}
As was the case in the power network, the node equations \eqref{eq:4} (with exception of the stochastic nodes) have no time dependency and can therefore be evaluated at each time step \(j\Delta t\).
Once again therefore we must evaluate the boundary conditions at each time step.

\subsection{Gas-power discretization}
\label{sec:gas-power-discr}
Again no further challenges arise in the discretization in the gas-power conversion plant equations \eqref{eq:5}.

\subsection{Stochastic process discretization}%
\label{sec:stoch-proc-discr}

The Ornstein-Uhlenbeck process is discretized with the explicit Euler-Maruyama
method, see~\cite{Saitostabilitystochastic1996}. Due to the explicit nature, the time steps for this method
must usually be chosen much finer than the time steps for the implicit box
scheme, as we detail below.  To make this distinction explicit we call the stepsize for the method
\(\Delta t_\text{stoch}\). To choose the boundary condition at time \(t^{*} = t
+\Delta t\), we make steps of size \(\Delta t_\text{stoch}\) according to
\begin{equation}
  \label{eq:9}
  P(t+\Delta t_\text{stoch})  = P(t)+ \theta(\hat{P}(t)-P(t))\Delta t_\text{stoch} + \sigma S(\Delta t_\text{stoch}),
\end{equation}
where \(\hat{P}\) takes on the role of the deterministic mean \(\mu\) of equation \eqref{eq:oup} and \(S(p)\) is a sample from a normal distribution with mean \(0\) and variance \(\Delta t_\text{stoch}\).
The same process is applied to get the discretized values of \(Q(t)\).
For stability in the mean (see again~\cite{Saitostabilitystochastic1996}), this discretization has the stepsize constraint
\begin{equation}
  \label{eq:10}
  \abs{1-\theta \Delta t_\text{stoch}} < 1,
\end{equation}
which for \(\theta>0\) yields
\begin{equation}
  \label{eq:11}
  0<\Delta t_\text{stoch} <  \frac{2}{\theta}\ .
\end{equation}

\begin{figure}[ht]
  \centering
  \includegraphics[width=0.85\textwidth]{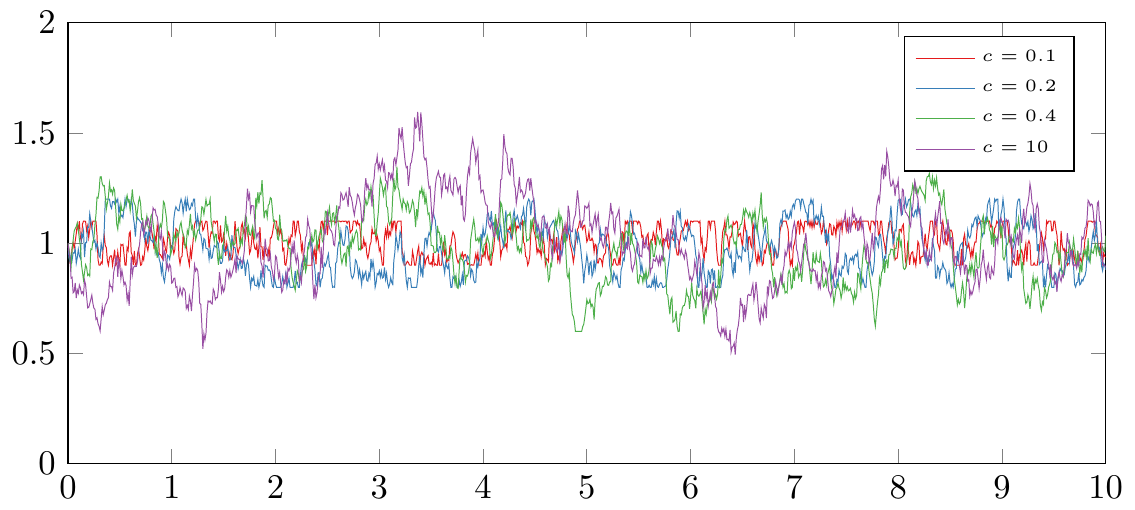}
  \caption{Ornstein-Uhlenbeck realizations for \(\mu = 1.0,\theta = 3.0,\sigma=0.45\) and different cut-off values \(c\).}
  \label{fig:cutoff-oup}
\end{figure}

\begin{figure}[ht]
  \centering
  \includegraphics[width=0.85\textwidth]{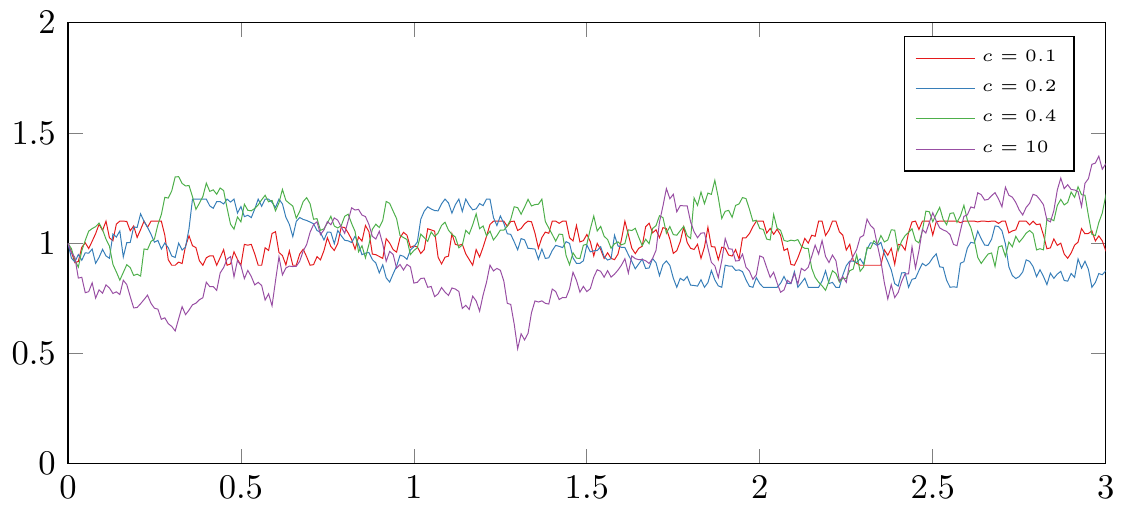}
  \caption{Zoomed-in version of \autoref{fig:cutoff-oup}.}
  \label{fig:cutoff-oup-zoom}
\end{figure}

In addition we also restrict the stochastic power demand according to
\begin{equation}
	\label{eq:7}
	\begin{aligned}
		(1-c) \hat{P}(t) \leq P(t) \leq (1+c) \hat{P}(t)\text{ if }\hat{P}(t)>0\\
		(1-c) \hat{P}(t) \geq P(t) \geq (1+c) \hat{P}(t)\text{ if }\hat{P}(t)<0
	\end{aligned}
\end{equation}
for some cutoff \(c\) with \(0\leq c \leq 1\). If the condition is violated, \(P(t)\) is
set to the boundary of the allowed interval. This cut-off prevents too great
outliers that are probably unrealistic and in addition prevent our numerical
methods from converging.
It may be argued that a stochastic process, whose samples must sometimes be cast away to yield usable solutions is a bad fit for its purpose.  Unfortunately we are not aware of a process that has been shown to be especially accurate for power fluctuations. However, an alternative might be the Jacobi process, as recently proposed in \cite{CoskunKorn2021}, which stays within a pre-defined interval.
Samples of the OUP for a couple of choices for the cut-off can be found in \autoref{fig:cutoff-oup} and a zoomed in version in \autoref{fig:cutoff-oup-zoom}.
In these figures the influence of the cut-off is easily seen.
The process for \(Q(t)\) is the same.

\newpage
\section{Software tool and computational results}
\label{sec:scenario}

\subsection{Software tool}
For the computations we use the network simulation tool \emph{grazer}\footnote{\url{https://github.com/eike-fokken/grazer}} written in C++.
It is an open source software suite developed at the Chair of Scientific Computing at the University of Mannheim. %
For the purpose of longterm usability the following design goals have been chosen:
\begin{itemize}
\item Easy installation
\item Full C++17-standard compliance with tested support for compilers GCC-9+, Clang-9+ and Visual Studio 2019+ (Other compilers are probably easy to use because of the standard compliance.)
\item Few external dependencies
\item High test coverage
\item Clean warning profile
\item Open Source License (AGPL 3.0).
\end{itemize}
The dependencies we have are Eigen, see \cite{eigenweb}, N. Lohmanns json library\footnote{\url{https://github.com/nlohmann/json}}, googletest\footnote{\url{https://github.com/google/googletest}}, pcg-random, see \cite{oneill:pcg2014} and CLI11\footnote{\url{https://github.com/CLIUtils/CLI11}}.

\subsubsection{Installation}
\label{sec:installation}
In order to build grazer you need three pieces of software: CMake\footnote{\url{https://cmake.org/}}, Git\footnote{\url{https://git-scm.com/}} and a C++17 capable C++ compiler, e.g. clang\footnote{\url{https://clang.llvm.org/}}, gcc\footnote{\url{https://gcc.gnu.org/}} or msvc\footnote{\url{https://visualstudio.microsoft.com/vs/}}.
Installation can be done by executing
\begin{lstlisting}[language=bash,caption={Installation},basicstyle=\scriptsize\ttfamily\color{blue}]
git clone https://github.com/eike-fokken/grazer.git
cd grazer
git submodule update --init --recursive --depth=1
cmake -DCMAKE_BUILD_TYPE=Release -S . -B release
\end{lstlisting}
Afterwards there is a grazer binary in .../grazer/release/src/Grazer called \emph{grazer} or \emph{grazer.exe} (on Windows).

\subsubsection{Usage}
\label{sec:usage}
Up to now grazer is a command-line application usable from any shell convenient, that is controlled by a number of input json files. In the medium term future it is planned to also support a python interface.

Grazer is used by pointing it to a directory with input json files.
\begin{lstlisting}[language=bash,caption={Calling grazer},basicstyle=\scriptsize\ttfamily\color{blue}]
grazer run data/one_pipeline
\end{lstlisting}
for example will run the problem defined in the directory \emph{data/one\textunderscore pipeline}.
The problem directory contains a subdirectory \emph{problem}, which holds the json files \emph{problem\textunderscore data.json, topology.json, boundary.json, initial.json} and \emph{control.json}.
Note that the layout of \emph{topology.json} was heavily inspired by the layout of GasLib files, see \cite{Gaslib}.

After solving the problem, an output file will be generated in \emph{data/one\textunderscore pipeline/output}.
This again a json file, so it can be read with almost all software.
For ease of use, some helper programs, compiled alongside grazer can be found in \emph{release/helper\textunderscore functions/}.
For example calling
\begin{lstlisting}[language=bash,caption={Calling grazer},basicstyle=\scriptsize\ttfamily\color{blue}]
generate_printing_csv data/one_pipeline/output/<outputfilename>  p_br1
\end{lstlisting}
will extract the json data into a csv file for usage with plotting tools. Helpers that import these into native formats of python are planned.

In addition json schemas can be generated and inserted into the jsons (Up to now with the exception of \emph{problem\textunderscore data.json}) with
\begin{lstlisting}[language=bash,caption={Calling grazer},basicstyle=\scriptsize\ttfamily\color{blue}]
grazer schema make-full-factory data/one_pipeline
grazer schema insert-key data/one_pipeline 
\end{lstlisting}
This has the advantage that json-aware editors help the users to only write jsons that are valid inputs for grazer which cuts down on bug searches.

As a final note on the usage, be aware that although grazer runs only sequentially, the output filename is chosen ``atomically'', meaning that two instances of grazer running in parallel will not interfere with each others output.
This is especially useful when executing many runs of stochastic problems in a Monte-Carlo method as was done in the present work.

Note that parts of the API are still subject to change.
For an up-to-date explanation check out the userguide in \emph{docs/userguide.tex} in the repository.

\subsubsection{A rough overview of the inner workings}
\label{sec:rough-overview-inner}
On execution grazer will read those files, configure the Newton solver according to settings in \emph{problem\textunderscore data.json}, construct a representation of the network from \emph{topology.json}, set initial and boundary values from the respective files and then start solving the problem time step per time step.
In each time step the model equations and their derivatives described above are evaluated to find a solution of them with Newton's method.
If successful, the solution is saved and the next time step is started.
If no solution can be found, the user is notified and all data computed in prior time steps is written to the output files.
If all time steps can be solved, all data is written out.

If a stochastic component is present in the network, a pseudo random number generator must be initialized with a seed. These are generated automatically or taken from the \emph{boundary.json} file, if a seed is present in there.

\subsubsection{Optimization}
\label{sec:optimization}
Grazer is also capable of computing optimal controls. To this end, certain components in a network can supply cost and constraint functions as well as their first derivatives. This information together with derivatives of the model equations with respect to states and controls is transformed via the adjoint method (see e.g. \cite{DissKolb} for an explanation) into derivatives with respect to the controls only.
The latter are then handed over by grazer to IPOPT (see \cite{Waechter2006}), that actually computes the optimal controls.  In our trials we have used the linear solver MUMPS (see \cite{MUMPS:1,MUMPS:2}), yet any other solver that can be interfaced with IPOPT could be used.  Note that no second derivatives are provided by grazer, which means that for the optimization only quasi-Newton methods are available.  A short example of this optimization is provided at the end of this work.
It is probably noteworthy that grazer is capable of handling constraint and control discretizations, that are coarser than the state discretizations detailed in \autoref{sec:discretization}. This is done by evaluating constraints only every \(n\)th time step, where \(n\) can be chosen by the user. Of course this can lead to constraint violations in between and therefore any such solution should be checked for such an occurrence.
The coarser control discretization is instead handled by interpolating controls linearly between two control discretization points.

\subsection{Scenario description}
\label{sec:scenario-1}
Here we describe the considered scenario of a combined power and gas network.
All data can be found in the git repo \url{https://github.com/eike-fokken/efficient_network-data.git}.

\subsection{Specification of the power network}
\label{sec:power-network}
\begin{figure}[p]
  \centering
  \includegraphics[width=\textwidth]{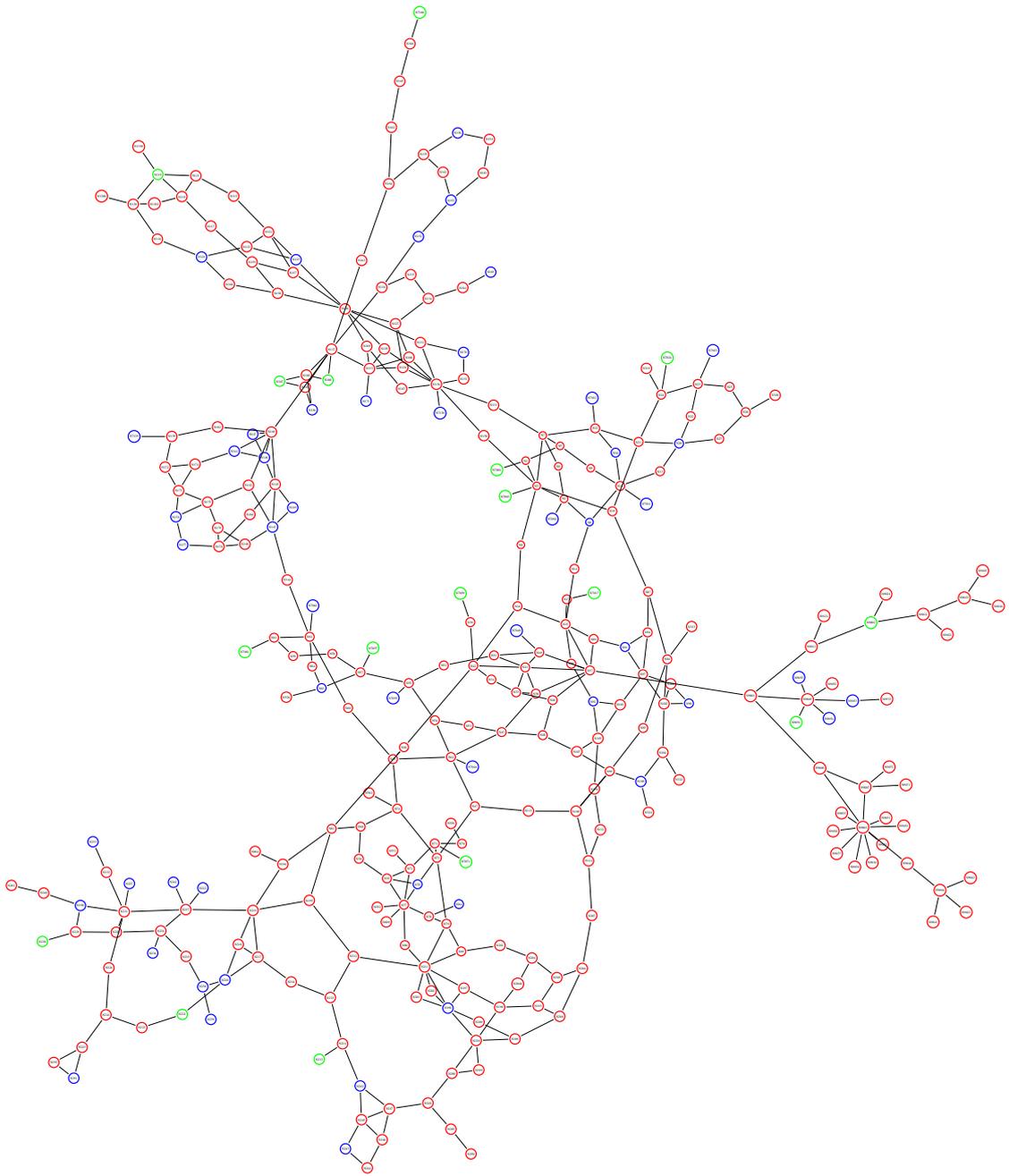}
  \caption{Power network with green gas powerplants, blue non-gas powerplants and red loads.}
  \label{fig:powernet}
\end{figure}

As starting point
for the power network we use the \emph{ieee-300-bus} system, as given in the
Matpowercase (see \cite{MATPOWER}) \emph{case300}. It is a power network of
\(69\) generator nodes, \(231\) load nodes and \(411\) transmission lines.
We alter the ieee-300 network in the following way.
\begin{itemize}
\item The power demand (real and reactive) is lowered by \(10\%\).
\item The former slack node N\(7049\) is changed into a PV-node.
\item The old PV-nodes given in table \autoref{tab:g2p} are turned into slack busses (V$\phi$-nodes).
\item All PQ-nodes are turned into stochastic PQ-nodes described in \autoref{sec:stoch-power-demand} and \autoref{sec:stoch-proc-discr}.
\end{itemize}
At the new V$\phi$-nodes, power that is generated from gas burned in gas power
plants is injected into the power network according to equation \eqref{eq:5}.
A picture of our power network can be found in \autoref{fig:powernet}.

\subsection{Specification of the gas network}%
\begin{figure}[h]
  \centering
  \includegraphics[width=\textwidth]{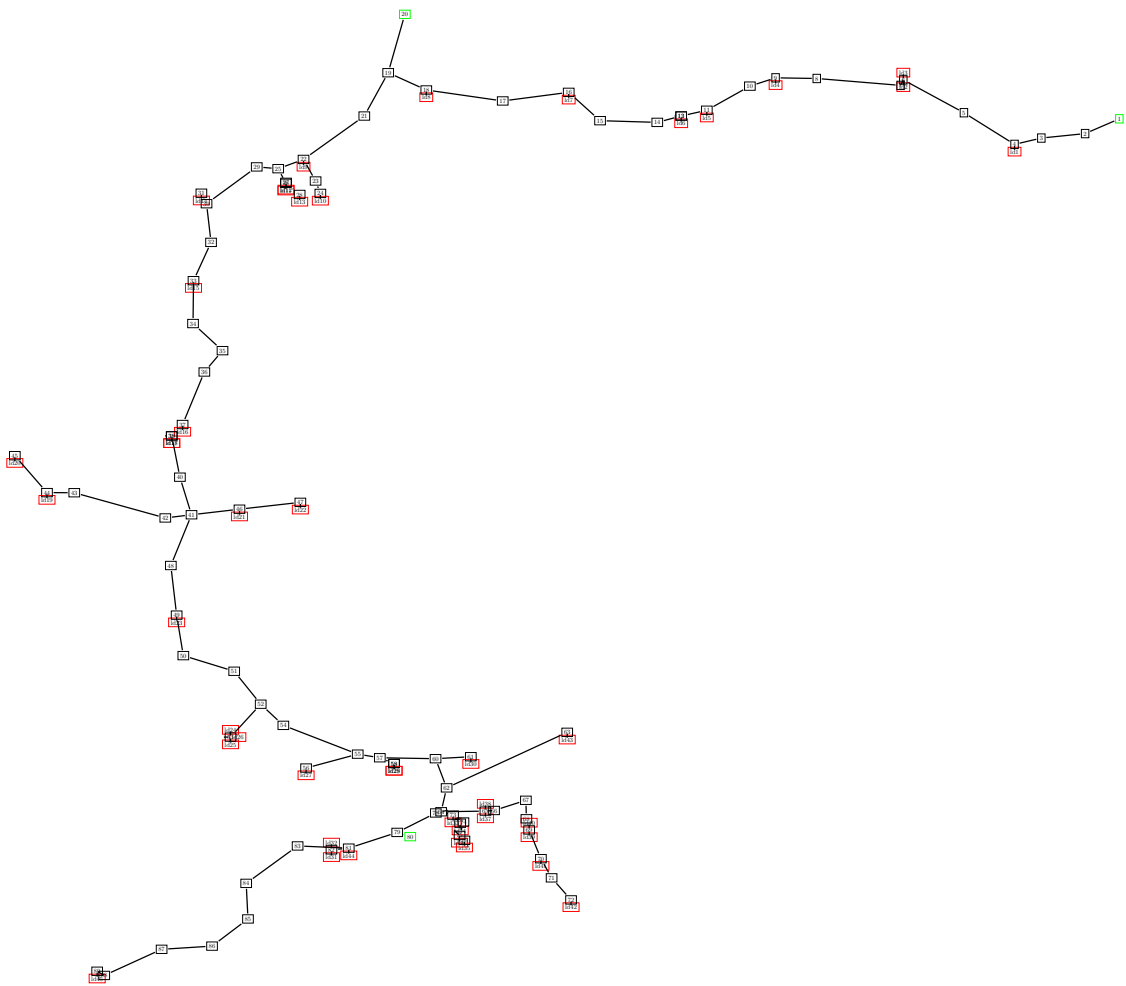}
  \caption{The gas network with green sources, red sinks and black junctions.}%
  \label{fig:gasnet}
\end{figure}

\label{sec:gas-network}
As starting point for the gas network we use the \emph{GasLib-134} system (see~\cite{Gaslib}).
A picture an be found in \autoref{fig:gasnet}.
\begin{table}[ht]
  \centering
  \begin{tabular}{lS[table-format=1.4,round-mode=places,round-precision=4]}
    \toprule
    Source node id& {inflow [\(\si{\meter \cubed \per \second}\)]}\\
    \midrule
    1 & 105.32815527751042\\
    20 & 280.6651734039139\\
    80 &    170.46067131857555\\
    \bottomrule
  \end{tabular}
  \caption{Inflow into the gas network}
  \label{tab:inflow_gas}
\end{table}

It is a gas network of \(86\) pipelines, \(3\) inflow nodes (\emph{sources}) and \(45\) outflow nodes (\emph{sinks}).
The inflow of gas remains constant over time and is given in Table \autoref{tab:inflow_gas}.

Here, 17 of the sinks, all gathered in table \autoref{tab:g2p},  draw gas to be converted into power.
The amount is set by the power network, and is computed from the power flow equations.  All other sinks do not consume gas.


\subsection{Specification of the Gas-Power connections}
The two networks are connected through gas-power conversion
plants, that turn gas into power, when power is needed and power into gas, when
surplus power is available. The gas-power conversion plants are arcs between the nodes listed in \autoref{tab:g2p}.
For simplicity they all share the same efficiencies both for power generation and gas generation, namely they have
\begin{equation}
  \label{eq:6}
  \begin{aligned}
    E_{\text{PtG}} &= \SI{43.56729}{\mega \watt \second \per \meter \cubed} \\
    E_{\text{GtP}} &= \SI{12.56}{\mega \watt \second \per \meter \cubed}\ .
  \end{aligned}
\end{equation}
For the smoothing constant we choose \(\kappa = \SI{60}{\meter \cubed \per \second}\).
As \(\kappa\) has no role but to mollify the kink in the switching from one conversion to the other,
the choice is purely driven by numerical factors.
The choice of \(\kappa\) must depend on the time step size,
where smaller time steps allow for smaller \(\kappa\) and therefore more sudden switching behavior.

All further data concerning these plants is gathered in
\autoref{tab:g2p}. There a real power demand is also given, which corresponds to
the default demand in our setting, when no uncertainty is present.
\begin{table}[ht]
  \centering
  \begin{tabular}[ht]{llS[table-format=1.4,round-mode=places,round-precision=4]}
    \toprule
    Gasnode & Powernode & { \(P[\SI{100}{\mega \watt}]\)} \\
    \midrule
    ld2  &N7017 &2.2890022500006006   \\
    ld6  &N7057 &1.3952492383902626   \\
    ld10 &N7071 &0.7217087852723401   \\
    ld12 &N7024 &2.7771189570764374   \\
    ld13 &N230  &2.5978326884895506   \\
    ld23 &N119  &19.299999999999976   \\
    ld24 &N221  &-0.08926590504578025 \\
    ld27 &N187  &11.402000000000005   \\
    ld29 &N7061 &2.7268692727170016   \\
    ld31 &N213  &2.0176362676661115   \\
    ld33 &N9051 &-0.3581000000000237  \\
    ld35 &N186  &11.402000000000001   \\
    ld36 &N7001 &2.1409910186666807   \\
    ld37 &N9002 &-0.0420000000000123  \\
    ld38 &N7166 &5.530000000000015    \\
    ld39 &N7003 &12.100000000000094   \\
    ld42 &N7039 &4.670240543306852    \\
    \bottomrule
  \end{tabular}
  \caption{Start and end nodes of gas-power-conversion plants and the deterministic demands of real power.}
  \label{tab:g2p}
\end{table}

\subsection{Specification of the stochastic power demand nodes}
\label{sec:spec-stoch-power}
As mentioned all PQ-nodes of the ieee-300-bus system are replaced by stochastic PQnodes.
As mean function we choose the power demands given by the ieee-300-bus problem, but lowered by \(\SI{10}{\percent}\).
In addition we choose for the drift coefficient \(\theta = 3 \), for the stability constraint we choose
\begin{equation}
  \label{eq:12}
  \Delta t_\text{stoch} =  \frac{0.1}{\theta},
\end{equation}
which results in rather high numbers of stochastic time steps but is unfortunately needed for convergence.
For the cutoff we use \(c = 0.4\).  The diffusion coefficient \(\sigma\) will be varied to compare different values.

\subsection{Computational results}
\label{sec:results}
In each run we simulate the combined network over the course of \SI{24}{\hour} (\SI{86400}{\second}) with a time stepsize of \SI{0.5}{\hour} (\SI{1800}{\second}).
\subsubsection{Steady state vs.\ stochastic example}
At first we simulate the network in a deterministic setting, which can be achieved by setting \(\sigma\) in \eqref{eq:9} to zero.
To keep the scenario simple, we choose steady-state initial conditions, which were generated by using arbitrary initial conditions and integrating them for a long time.  The resulting end state is then used as initial conditions for our setting.

In this deterministic setting we find the (constant) power demands in the gas plants given in \autoref{tab:g2p}.
To illustrate our results we will usually picture the situation of pipe \emph{p\textunderscore{}br71}, which is located to the lower right in \autoref{fig:gasnet} connecting nodes \(71\) and \(72\).
The steady-state solution in pipe \emph{p\textunderscore{}br71} remains constant over time as is fitting for a steady state solution.
The same is true for the flow, and also for all other pipes in the network.

Along with the deterministic setting, we simulate a scenario with \(\theta = 3.0\) and \( \sigma=0.45\) for all PQ-nodes.
The number of stochastic steps is set to at least \(1000\) which, due to stability constraints mentioned in \eqref{eq:11}, was then automatically raised to \(18000\).

\begin{figure}[ht]
  \centering
  \includegraphics[width=0.8\linewidth]{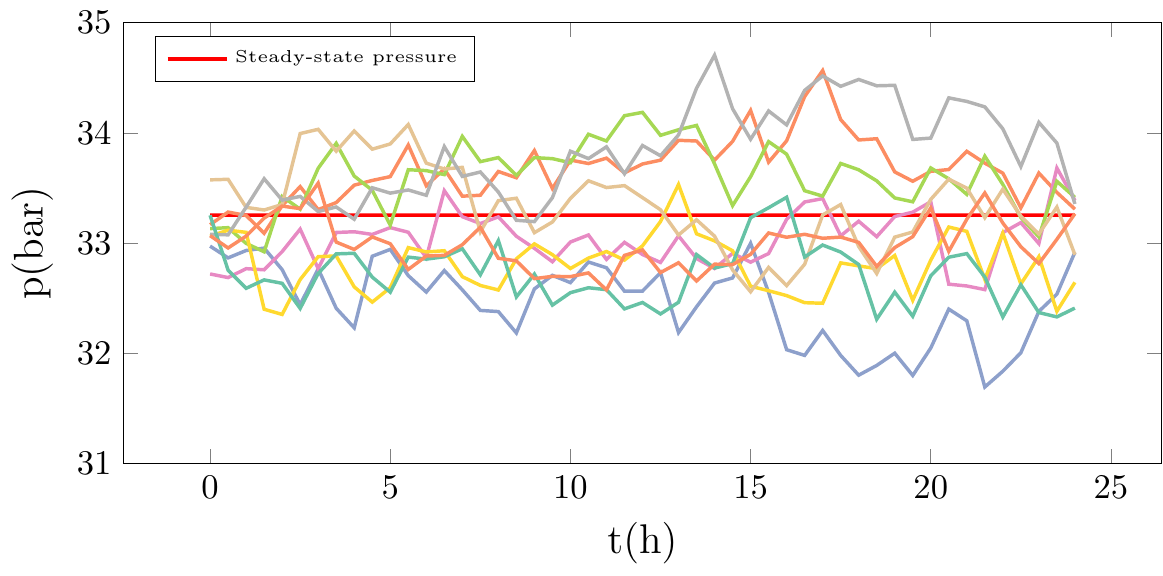}
  \caption{Pressure evolution in \emph{p\textunderscore{}br71} for deterministic and some realizations of stochastic demand.}
  \label{fig:singlepressure}
\end{figure}

\begin{figure}[ht]
  \centering
  \includegraphics[width=0.8\linewidth]{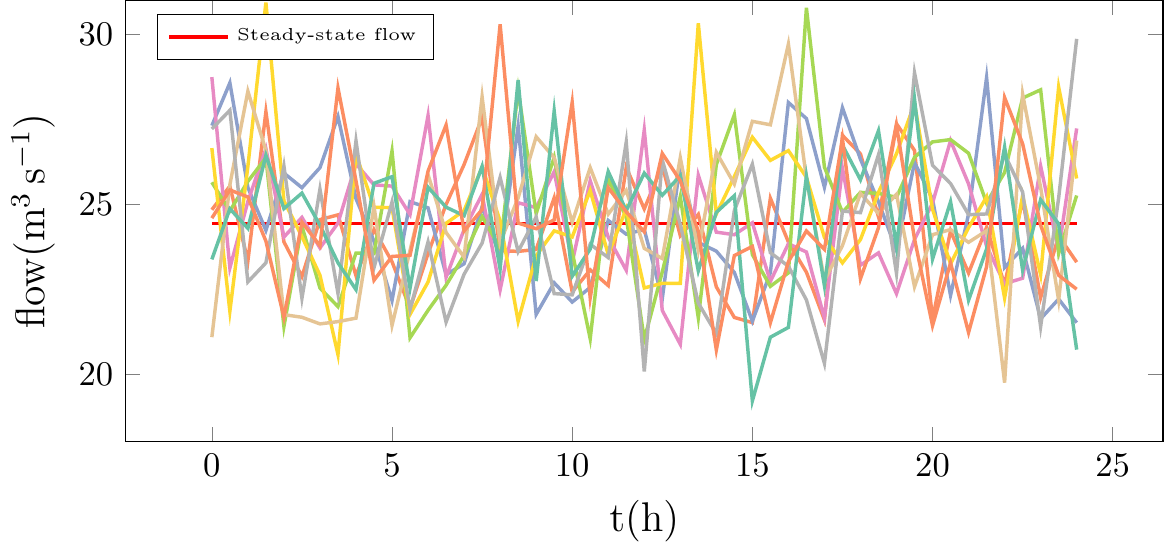}
  \caption{Flow evolution in \emph{p\textunderscore{}br71} for deterministic and some realizations of stochastic demand.}
  \label{fig:singleflow}
\end{figure}

A comparison of the steady-state and stochastic pressure can be seen in \autoref{fig:singlepressure} while
a comparison of the fluxes is given in \autoref{fig:singleflow}.


\begin{figure}[h!]
  \centering
  \includegraphics[width=0.8\linewidth]{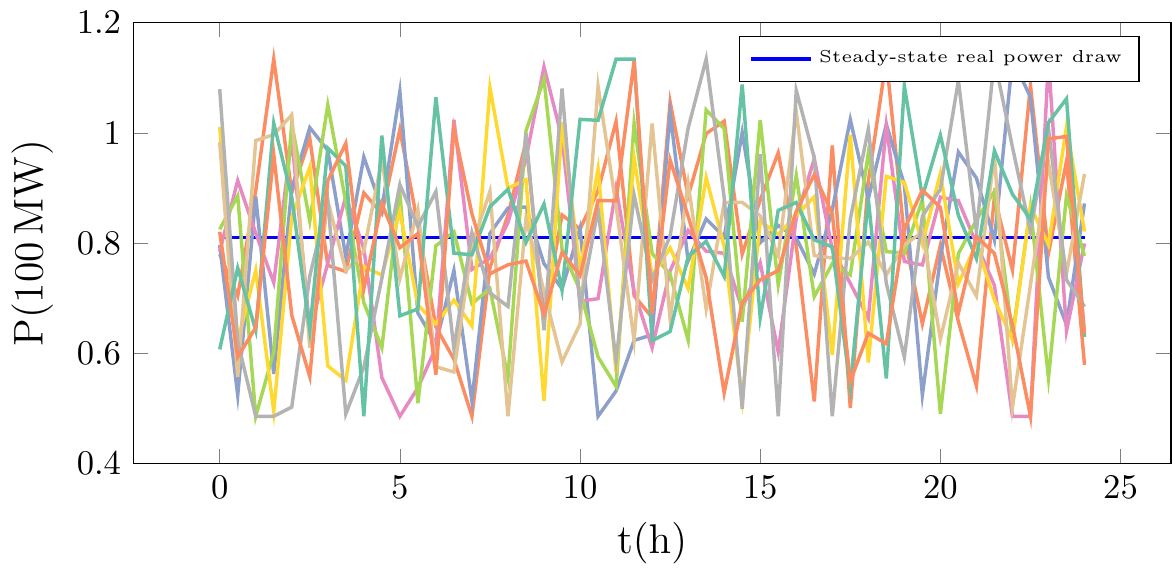}
  \caption{Real power demand in \emph{N1} for deterministic and stochastic demand with \(\sigma =0.45\).}
  \label{fig:realpower}
\end{figure}

\begin{figure}[h!]
  \centering
  \includegraphics[width=0.8\linewidth]{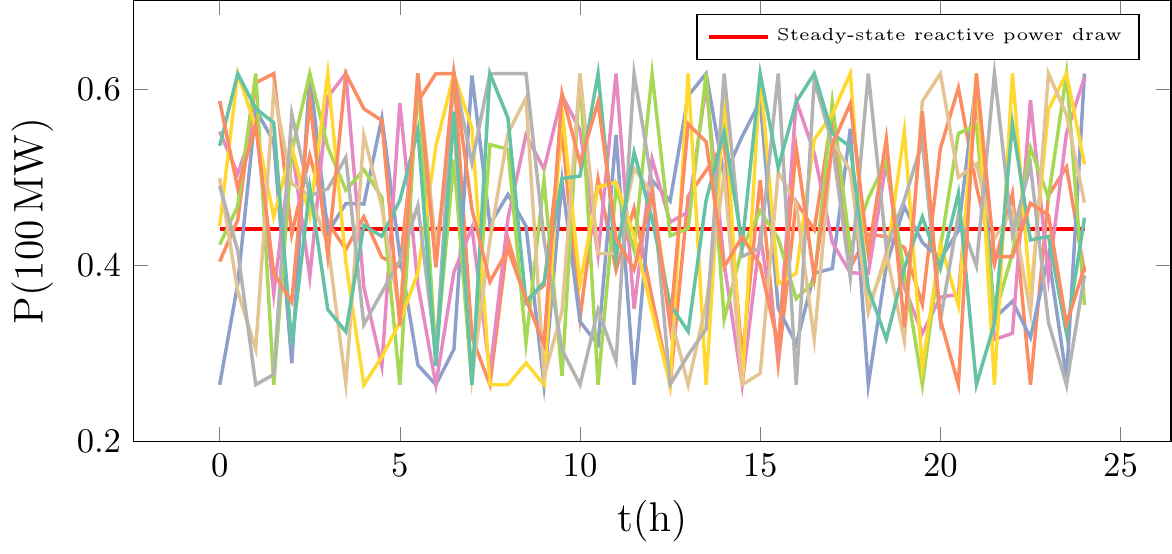}
  \caption{Reactive power demand in \emph{N1} for deterministic and stochastic demand with \(\sigma =0.45\).}
  \label{fig:reactivepower}
\end{figure}

In the power network we find for the PQ-node \emph{N1} power demands over time like those in \autoref{fig:realpower} and \autoref{fig:reactivepower}.  Of course the situation is similar for all PQ-nodes.

\subsubsection{Stochastic demand with variable noise }
\begin{figure}[ht]
  \centering
  \includegraphics[width=0.8\linewidth]{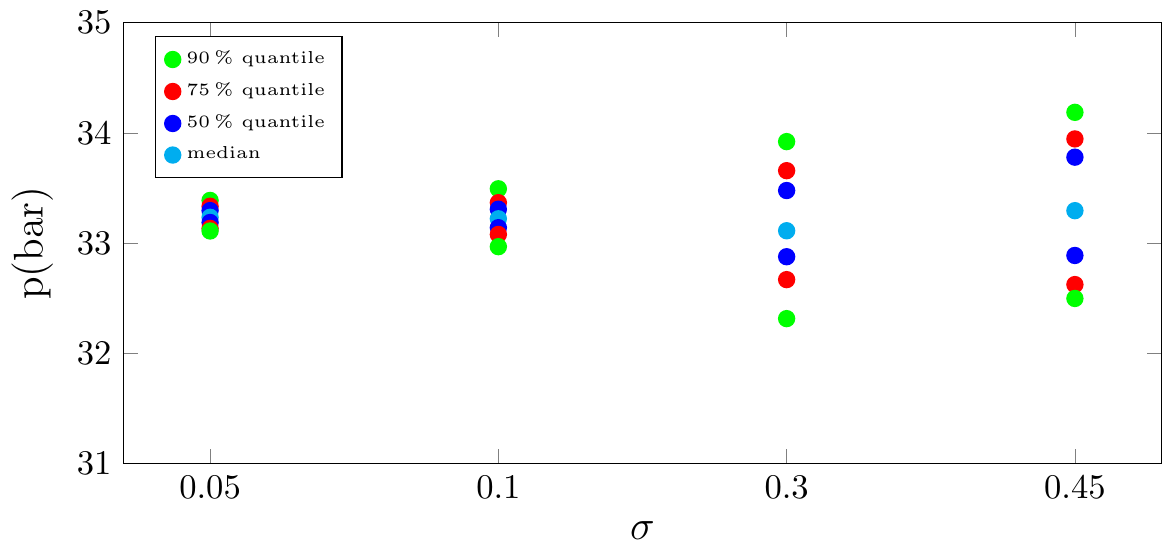}
  \caption{Comparison of pressure quantile boundaries at different \(\sigma\) at \(t = \SI{12}{\hour}\) in pipeline \emph{p\textunderscore{}br71}.}
  \label{fig:pressure_quantiles_12hours}
\end{figure}

\begin{figure}[ht]
  \centering
  \includegraphics[width=0.8\linewidth]{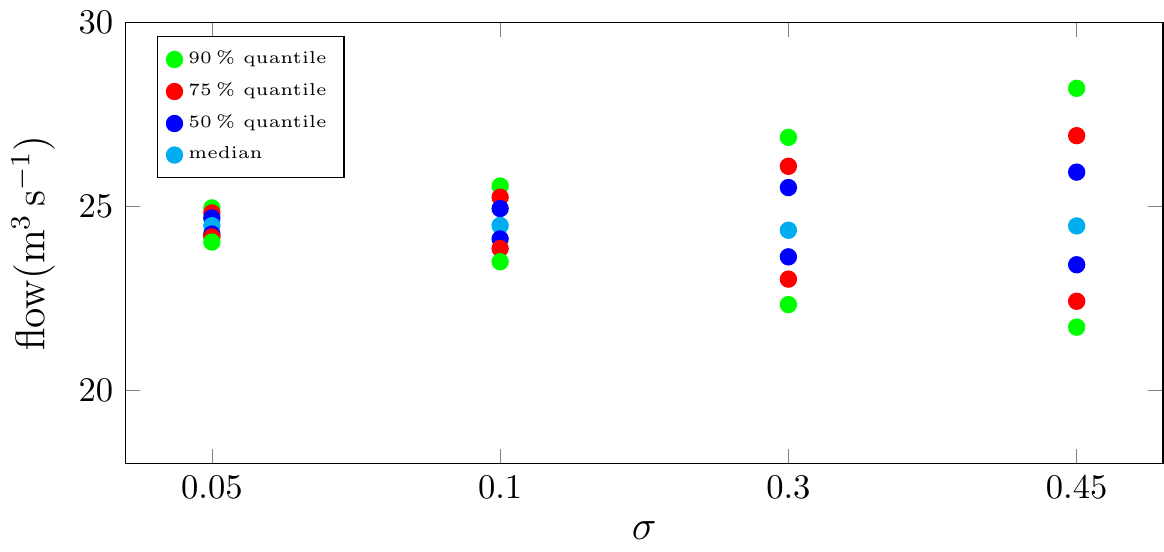}
  \caption{Comparison of flow quantile boundaries at different \(\sigma\) at \(t = \SI{12}{\hour}\) in pipeline \emph{p\textunderscore{}br71}.}
  \label{fig:flow_quantiles_12hours}
\end{figure}

Now we examine repercussions of the uncertainty on the gas network.
Therefore we make \(100\) runs for each \(\sigma \in \{0.05,0.1,0.3,0.45\}\) and compare the quantiles at \SI{50}{\percent}, \SI{75}{\percent} and \SI{90}{\percent}.
Taking an arbitrary point in time, \(t = \SI{12}{\hour}\), the quantiles for the pressure can be seen in \autoref{fig:pressure_quantiles_12hours}.
For the flow the quantile comparison can be found in \autoref{fig:flow_quantiles_12hours}.

For both quantities we see the expected expansion of quantile boundaries with higher diffusion \(\sigma\).

\subsubsection{Comparison of deterministic and stochastic pressure prediction }
\begin{figure}[p]
  \centering
  \includegraphics[width=1.0\textwidth]{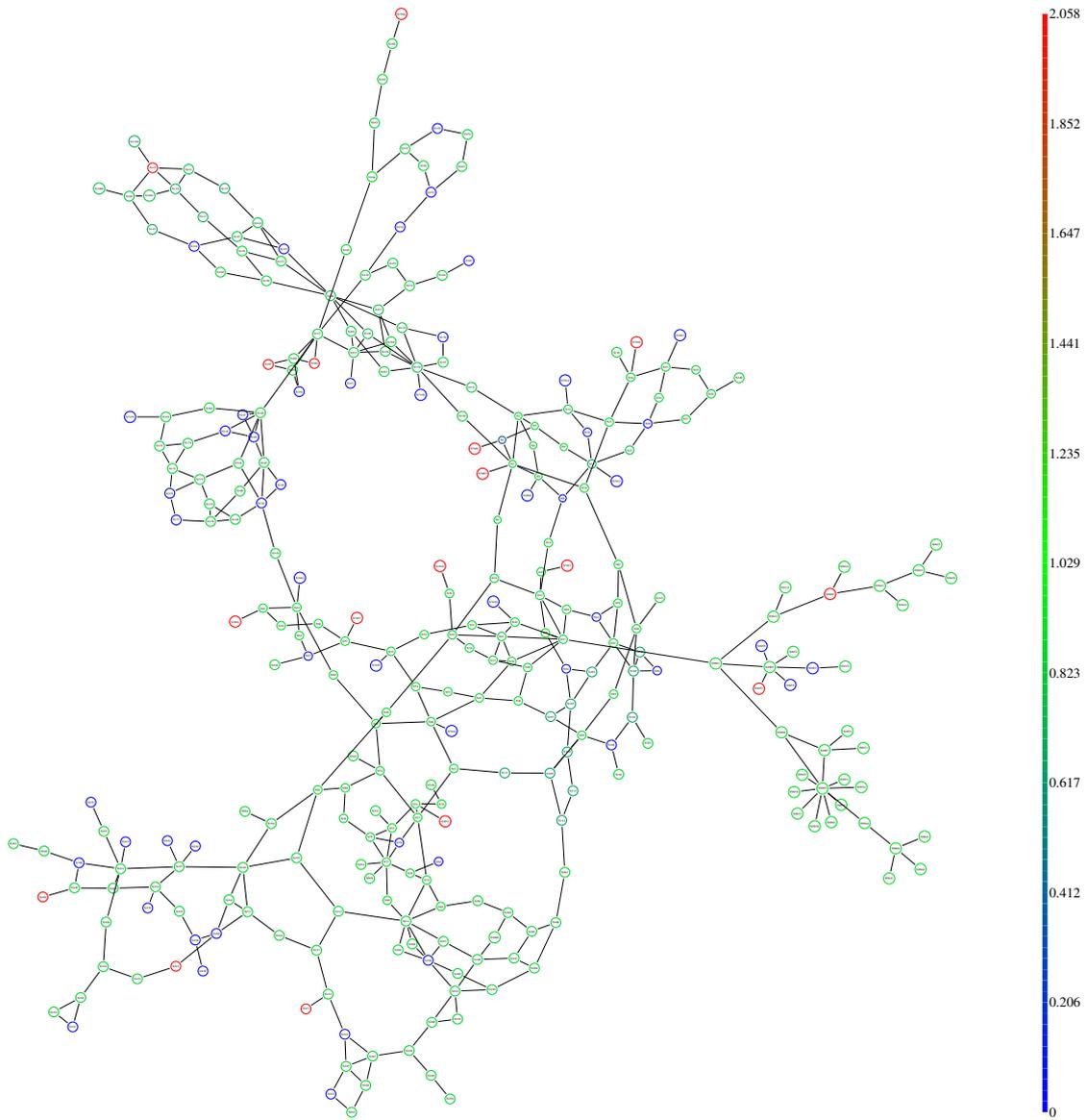}
  \caption{Heatmap of the maximal real power deviation over the course of \SI{24}{\hour}, units are in \SI{100}{\mega\watt}.}
  \label{fig:heatmap_power_difference}
\end{figure}

\begin{figure}[p]
  \centering
  \includegraphics[width=1.0\textwidth]{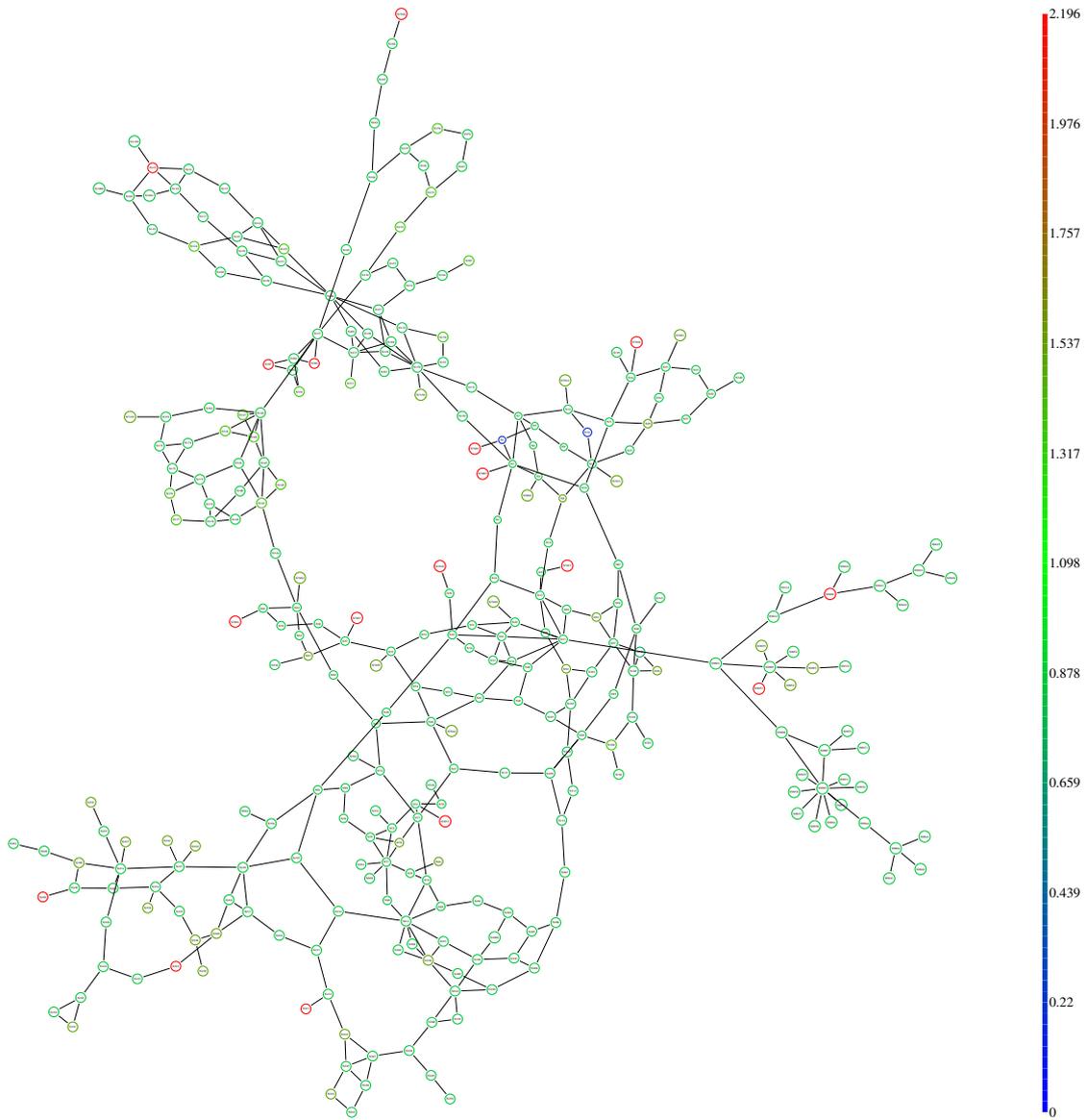}
  \caption{Heatmap of the maximal reactive power deviation over the course of \SI{24}{\hour}, units are in \SI{100}{\mega\watt}.}
  \label{fig:heatmap_reactive_power_difference}
\end{figure}

Now we give an overview of the impact of the volatility in power demand on the network.
Therefore we revisit the scenario with the highest volatility, that is with \(\sigma = 0.45\) and consider again a time frame of \SI{24}{\hour}.
In \autoref{fig:heatmap_power_difference} one can see the maximal deviation of real power demand from the steady state solution.
At first glance this looks similar to \autoref{fig:powernet}, just with colors cycled around.
This is due to the fact, that the load nodes have defined volatility as they follow their own Ornstein-Uhlenbeck process approximation defined in equation \eqref{eq:9}.
The PV-nodes on the other hand have zero volatility in real power.
Yet the V$\phi$-nodes must account for all remaining power demand and as such have the highest volatility.
A similar picture can be found in \autoref{fig:heatmap_reactive_power_difference}, where the deviation of the reactive power is depicted.
Here the PV-nodes do not have zero volatility, yet it seems that they also do not carry much volatility in \(Q\).

\begin{figure}[p]
  \centering
  \includegraphics[width=1.0\textwidth]{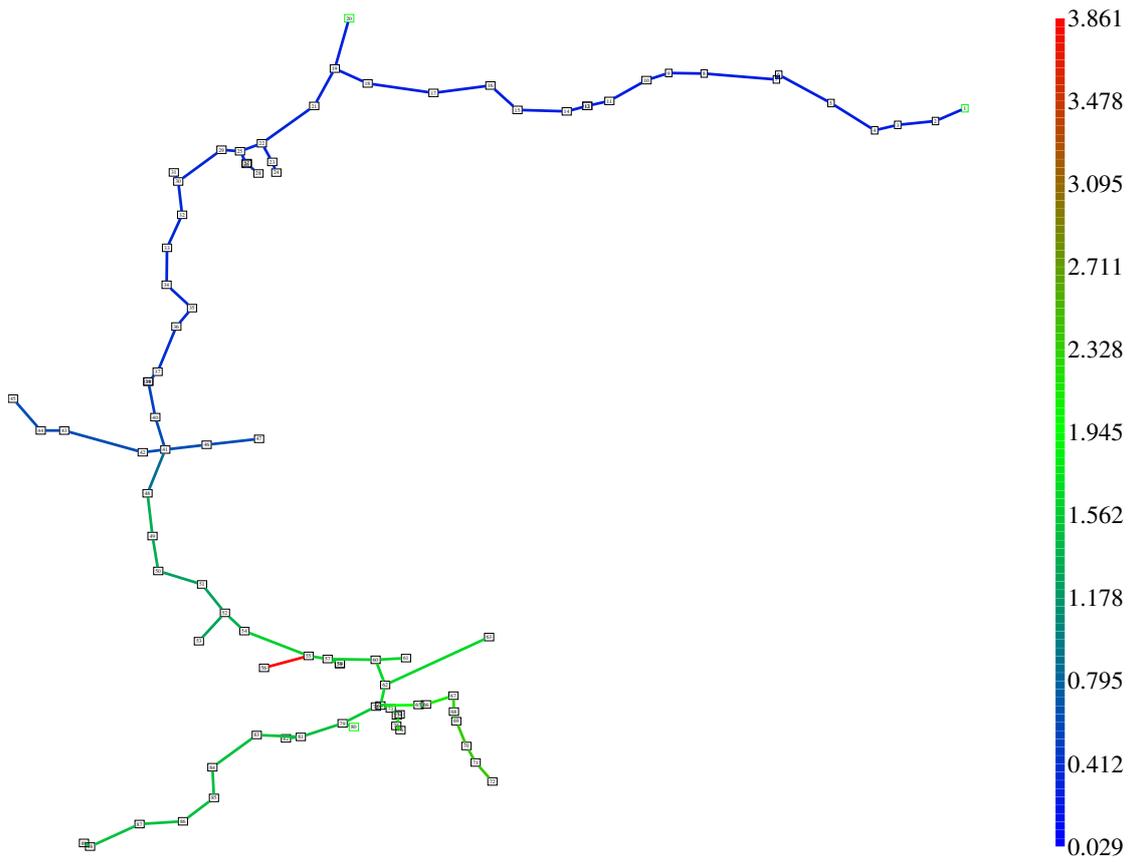}
  \caption{Heatmap of the maximal pressure deviation over the course of \SI{24}{\hour}, units are in \si{\bar}.}
  \label{fig:heatmap_pressure_difference}
\end{figure}

At last we consider the possible impact of the volatility in the power network on the gas network.
Therefore we show the maximal pressure deviation from the steady-state solution over the course of \(\SI{24}{\hour}\) for \(\sigma = 0.45\) in \autoref{fig:heatmap_pressure_difference}.
It is easily seen that the lower part of the network experiences much higher pressure volatility than the upper part.
This is expected, as on the one hand the upper part has higher pressure as the three gas sources are located there and on the other hand many more gas-power-conversion plants are located in the lower part, so that the volatility can add up.


\subsubsection{Optimization example}%
\label{sec:optimization-example}
Finally, we show results of an optimization task carried out with grazer.
The Gaslib-134 network actually contains two controllable components, a compressor between the nodes \(29\)/\(30\) and a control valve between the nodes \(65\)/\(66\).
We take the steady state solution from above but add two continuous constraints in order to make the controllable components actually do some work.
At the sink ld\textunderscore{}22 we impose a lower pressure bound of \(\SI{70}{\bar}\) at \(t=0\), \(\SI{90}{\bar}\) at \(t=\SI{24}{\hour}\) and interpolate linearly in between.
In addition we impose an upper pressure bound at sink ld\textunderscore{}40 of \(\SI{90}{\bar}\) at \(t=0\), \(\SI{70}{\bar}\) at \(t=\SI{24}{\hour}\) and again interpolate linearly in between.
As cost function we choose
\begin{eqnarray*}
  \int_0^{\SI{24}{\hour}}u_{\text{compressor}}(t)^2\dif t,
\end{eqnarray*}
so that using the valve is free but compressor costs should be minimized. The control is discretized with the same discretization already used by the states, yielding 49 timepoints.
Also the constraints are evaluated at every state timestep.
While grazer is capable of using coarser resolutions of both constraints and controls, the problem at hand is small enough to compute a solution in approximately a minute on a workstation.  Using only \(11\) controls cut this time in half.  In addition evaluating the constraints only every fifth step reduces the time again by half.
The control of the valve is constrained to not exceed \(\SI{40}{\bar}\) to keep the optimization routine from trying controls that are too high to yield a solution of the simulation.
The compressor control is capped at \(\SI{120}{\bar}\), although this bound is never attained.

\begin{figure}[h!]
  \centering
  \includegraphics[width=0.8\textwidth]{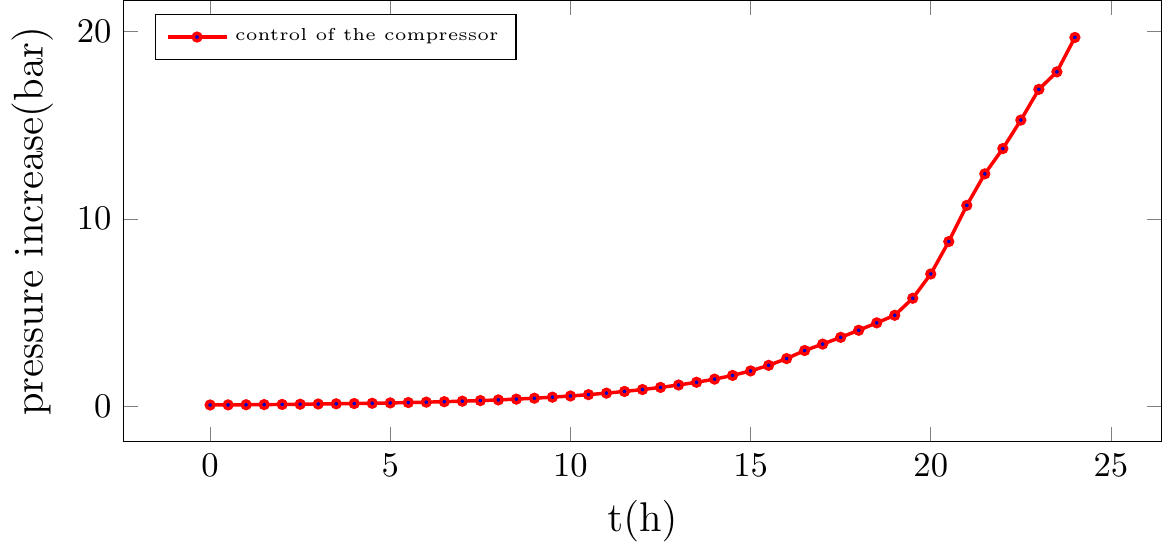}
  \caption{Computed optimal control of the compressor at nodes \(29\)/\(30\).}
  \label{fig:compressor_control}
\end{figure}

\begin{figure}[h!]
  \centering
  \includegraphics[width=0.8\textwidth]{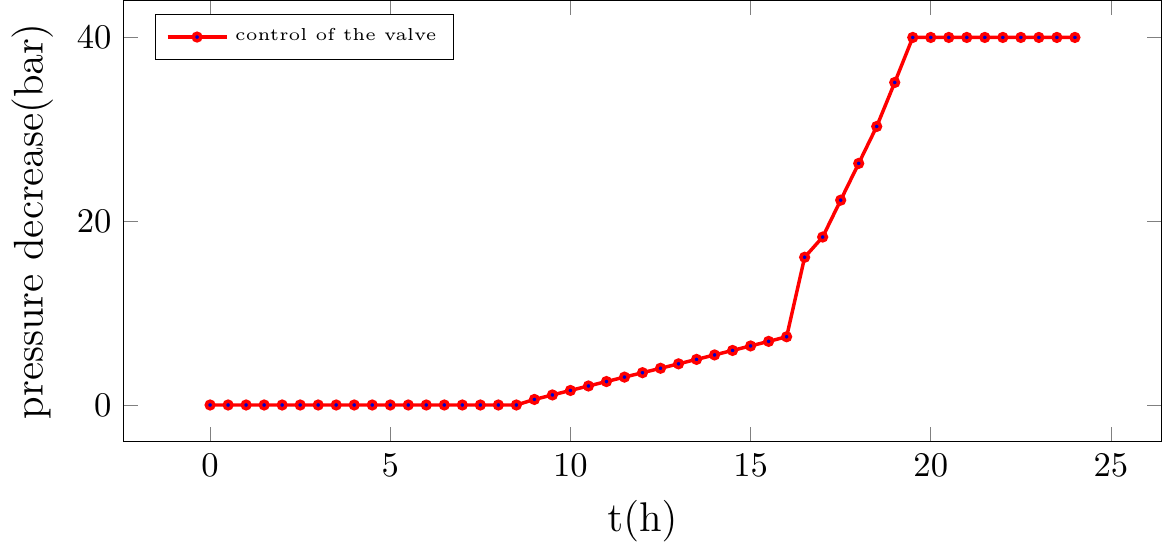}
  \caption{Computed optimal control of the valve at nodes \(65\)/\(66\).}
  \label{fig:valve_control}
\end{figure}

With this data grazer computes the optimal controls in \autoref{fig:compressor_control} and \autoref{fig:valve_control}.
The compressor control in \autoref{fig:compressor_control} nicely ramps up as the lower pressure bound in ld\textunderscore{}22 rises.
On the other hand, the valve control in \autoref{fig:valve_control} stays at zero until this is not sufficient anymore to satisfy the decreasing upper pressure bound in ld\textunderscore{}40 at which point the control rises up to the maximal value, staying there until the end. As the valve control incurs no cost, this is one of many possible configurations.

\begin{figure}[h!]
  \centering
  \includegraphics[width=0.8\textwidth]{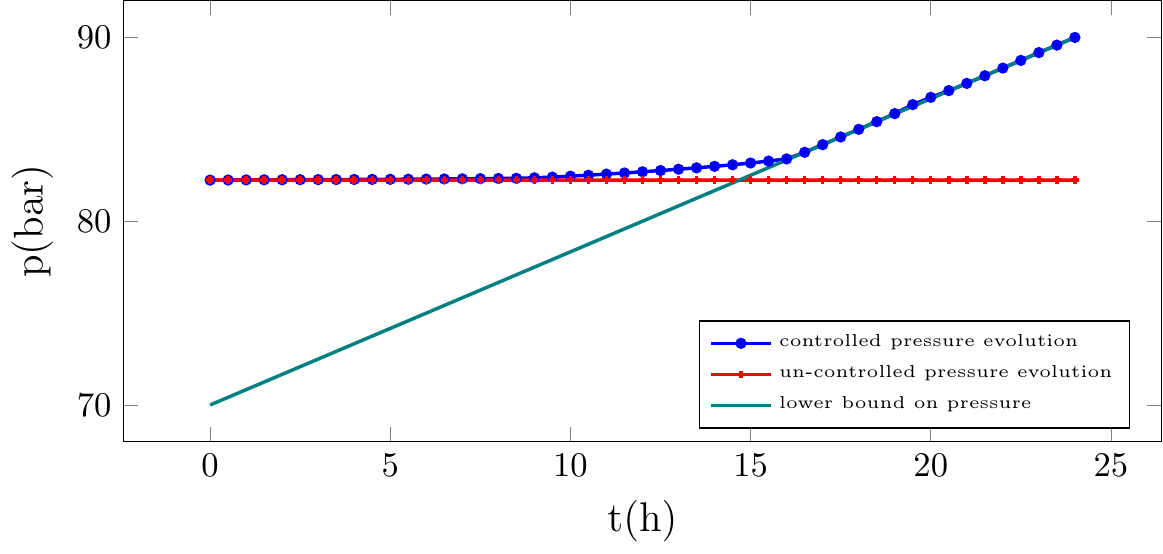}
  \caption{Comparison of controlled and uncontrolled pressure at ld\textunderscore{}22.}
  \label{fig:lowerbound}
\end{figure}

\begin{figure}[h!]
  \centering
  \includegraphics[width=0.8\textwidth]{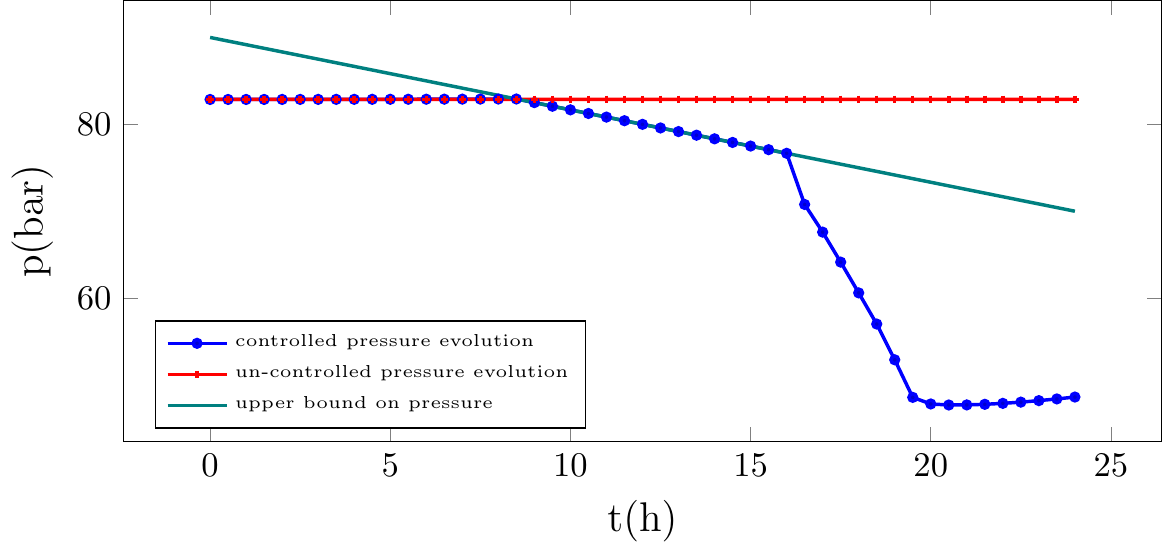}
  \caption{Comparison of controlled and uncontrolled pressure at ld\textunderscore{}40.}
  \label{fig:upperbound}
\end{figure}

A comparison of pressure evolution at the two sinks ld\textunderscore{}22 and ld\textunderscore{}40 is given in \autoref{fig:lowerbound} and \autoref{fig:upperbound}.
It can be seen that the compressor increases the pressure just enough to satisfy the lower pressure bound as its usage is penalized, while the (free to use) valve at first matches the upper pressure bound exactly but later on over-compensates it rather strongly.

\section{Summary and future work}
\label{sec:summary-future-work}
We introduced the new open-source software tool \emph{grazer} that can be used
to efficiently simulate numerical problems that are defined on networks. We used grazer to
simulate a coupled gas and power network with uncertain power demand presented
repercussions of the uncertain power demand within the gas network.

Future work includes the extension of grazer to more complex optimization problems and other
types of uncertainty. %

\subsection*{Acknowledgments}
The authors acknowledge support through the BMBF project ENets (05M18VMA) and 320021702/GRK2326, 333849990/IRTG-2379 and 
DFG HE5386/18-1 and DFG HE5386/23-1.

\newpage
\printbibliography[heading=bibintoc]

\end{document}